\documentclass{article}[12pt]
\usepackage[T1]{fontenc}
\usepackage[latin1]{inputenc}
\usepackage[english]{babel}
\usepackage{hyperref}
\usepackage[dvips]{graphics}
\usepackage{amsmath}
\usepackage{amssymb}
\usepackage{amsopn}
\usepackage{amsbsy}
\usepackage{amstext}
\usepackage{latexsym}
\usepackage{amsfonts}
\usepackage{array}
\usepackage{epsfig}
\usepackage{mathrsfs}
\title{New results on the least common multiple of consecutive integers}
\author{Bakir FARHI and Daniel KANE}
\date{}

\newtheorem{thm}{Theorem}[section]
\newtheorem{prop}[thm]{Proposition}
\newtheorem{lemma}[thm]{Lemma}

\newtheorem{coll}[thm]{Corollary}

\newcommand{\floor}[1]{\left\lfloor #1 \right\rfloor}

\let\epsilon=\varepsilon

\def\gcd{{\rm gcd}}
\def\lcm{{\rm lcm}}

\def\EMdash{\leavevmode\hbox to 7.5mm{\vrule height .63ex depth -.59ex
    width 5.4mm\hfill}}

\begin{document}
\maketitle
\begin{center}
Département de Mathématiques, Université du Maine, \\
Avenue Olivier Messiaen, 72085 Le Mans Cedex 9, France. \\
bakir.farhi@gmail.com
\end{center}
\begin{center}
Harvard University, Department of Mathematics \\
1 Oxford Street, Cambridge MA 02139, USA. \\
aladkeenin@gmail.com
\end{center}
\begin{abstract}
When studying the least common multiple of some finite sequences
of integers, the first author introduced the interesting
arithmetic functions $g_k$ $(k \in \mathbb{N})$, defined by
$g_k(n) := \frac{n (n + 1) \dots (n + k)}{\lcm(n , n + 1 , \dots ,
n + k)}$ $(\forall n \in \mathbb{N} \setminus \{0\})$. He proved
that $g_k$ $(k \in \mathbb{N})$ is periodic and $k!$ is a period
of $g_k$. He raised the open problem consisting to determine the
smallest positive period $P_k$ of $g_k$. Very recently, S. Hong
and Y. Yang have improved the period $k!$ of $g_k$ to $\lcm(1 , 2
, \dots , k)$. In addition, they have conjectured that $P_k$ is
always a multiple of the positive integer $\frac{\lcm(1 , 2 ,
\dots , k , k + 1)}{k + 1}$. An immediate consequence of this
conjecture states that if $(k + 1)$ is prime then the exact period
of $g_k$ is precisely equal to $\lcm(1 , 2 , \dots , k)$.

In this paper, we first prove the conjecture of S. Hong and Y.
Yang and then we give the exact value of $P_k$ $(k \in
\mathbb{N})$. We deduce, as a corollary, that $P_k$ is equal to
the part of $\lcm(1 , 2 , \dots , k)$ not divisible by some prime.
\end{abstract}
{\it MSC:} 11A05~\vspace{1mm}\\
{\it Keywords:} Least common multiple; arithmetic function;
exact period.
\section{Introduction}~\vspace{-1mm}

Throughout this paper, we let $\mathbb{N}^*$ denote the set
$\mathbb{N} \setminus \{0\}$ of positive integers.

Many results concerning the least common multiple of sequences of
integers are known. The most famous is nothing else than an
equivalent of the prime number theorem; it sates that $\log\lcm(1
, 2 , \dots , n) \sim n$ as $n$ tends to infinity (see e.g.,
\cite{hw}). Effective bounds for $\lcm(1 , 2 , \dots , n)$ are
also given by several authors (see e.g., \cite{ha} and \cite{n}).

Recently, the topic has undergone important developments. In
\cite{b}, Bateman, Kalb and Stenger have obtained an equivalent
for $\log\lcm(u_1 , u_2 , \dots , u_n)$ when ${(u_n)}_n$ is an
arithmetic progression. In \cite{c}, Cilleruelo has obtained a
simple equivalent for the least common multiple of a quadratic
progression. For the effective bounds, Farhi \cite{f1} \cite{f2}
got lower bounds for $\lcm(u_0 , u_1 , \dots , u_n)$ in both cases
when ${(u_n)}_n$ is an arithmetic progression or when it is a
quadratic progression. In the case of arithmetic progressions,
Hong and Feng \cite{hf} and Hong and Yang \cite{hy} obtained some
improvements of Farhi's lower bounds.

Among the arithmetic progressions, the sequences of consecutive
integers are the most well-known with regards the properties of
their least common multiple. In \cite{f2}, Farhi introduced the
arithmetic function $g_k : \mathbb{N}^* \rightarrow
\mathbb{N}^*$ $(k \in \mathbb{N})$ which is defined by:
$$g_k(n) := \frac{n (n + 1) \dots (n + k)}{\lcm(n , n + 1 , \dots , n + k)} ~~~~~~ (\forall n \in \mathbb{N}^*) .$$
Farhi proved that the sequence ${(g_k)}_{k \in \mathbb{N}}$
satisfies the recursive relation:
\begin{equation}\label{eqf}
g_k(n) = \gcd\left(k! , (n + k) g_{k - 1}(n)\right) ~~~~~~
(\forall k , n \in \mathbb{N}^*) .
\end{equation}
Then, using this relation, he deduced (by induction on k) that
$g_k$ $(k \in \mathbb{N})$ is periodic and $k!$ is a period of
$g_k$. A natural open problem raised in \cite{f2} consists to
determine the exact period (i.e., the smallest positive period) of
$g_k$.

For the following, let $P_k$ denote the exact period of $g_k$. So,
Farhi's result amounts that $P_k$ divides $k!$ for all $k \in
\mathbb{N}$. Very recently, Hong and Yang have shown that $P_k$
divides $\lcm(1 , 2 , \dots , k)$. This improves Farhi's result
but it doesn't solve the raised problem of determining the
$P_k$'s. In their paper \cite{hy}, Hong and Yang have also
conjectured that $P_k$ is a multiple of $\frac{\lcm(1 , 2 , \dots
, k + 1)}{k + 1}$ for all nonnegative integer $k$. According to
the property that $P_k$ divides $\lcm(1 , 2 , \dots , k)$
$(\forall k \in \mathbb{N})$, this conjecture implies that the
equality $P_k = \lcm(1 , 2 , \dots , k)$ holds at least when $(k +
1)$ is prime.

In this paper, we first prove the conjecture of Hong and Yang and
then we give the exact value of $P_k$ $(\forall k \in
\mathbb{N})$. As a corollary, we show that $P_k$ is equal to the
part of $\lcm(1 , 2 , \dots , k)$ not divisible by some prime and
that the equality $P_k = \lcm(1 , 2 , \dots , k)$ holds for an
infinitely many $k \in \mathbb{N}$ for which $(k + 1)$ is not
prime.

\section{Proof of the conjecture of Hong and Yang}~

We begin by extending the functions $g_k$ $(k \in \mathbb{N})$
to $\mathbb{Z}$ as follows:\\
$\bullet$ We define $g_0 : \mathbb{Z} \rightarrow \mathbb{N}^*$ by
$g_0(n) = 1$, $\forall n \in \mathbb{Z}$. \\
$\bullet$ If, for some $k \geq 1$, $g_{k - 1}$ is defined, then we
define $g_k$ by the relation:
\begin{equation}\label{eqf2}
g_k(n) = \gcd\left(k! , (n + k) g_{k - 1}(n)\right) ~~~~~~
(\forall n \in \mathbb{Z}) . \tag{$\ref{eqf}'$}
\end{equation}
Those extensions are easily seen to be periodic and to
have the same period as their restriction to $\mathbb{N}^*$. The
following proposition plays a vital role in what follows:
\begin{prop}\label{p}
For any $k \in \mathbb{N}$,we have $g_k(0) = k!$.
\end{prop}
{\bf Proof.} This follows by induction on $k$ with using the
relation
(\ref{eqf2}).\penalty-20\null\hfill$\blacksquare$\par\medbreak

We now arrive at the theorem implying the conjecture of Hong and
Yang.
\begin{thm}\label{t7}
For all $k \in \mathbb{N}$, we have:
$$P_k = \frac{\lcm(1 , 2 , \dots , k + 1)}{k + 1} . \gcd\left(P_k + k + 1 , \lcm(P_k + 1 , P_k + 2 , \dots ,
P_k + k)\right) .$$
\end{thm}
The proof of this theorem needs the following lemma:
\begin{lemma}\label{l1}
For all $k \in \mathbb{N}$, we have:
$$\lcm(P_k , P_k + 1 , \dots , P_k + k) = \lcm(P_k + 1 , P_k + 2 , \dots , P_k + k) .$$
\end{lemma}
{\bf Proof of the Lemma.} Let $k \in \mathbb{N}$ fixed. The
required equality of the lemma is clearly equivalent to say that
$P_k$ divides $\lcm(P_k + 1 , P_k + 2 , \dots , P_k + k)$. This
amounts to showing that for any prime number $p$:
\begin{equation}\label{eq2}
v_p(P_k) \leq v_p\left(\lcm(P_k + 1 , \dots , P_k + k)\right) =
\max_{1 \leq i \leq k} v_p(P_k + i).
\end{equation}
So it remains to show (\ref{eq2}). Let $p$ be a prime number.
Because $P_k$ divides $\lcm(1 , 2 , \dots , k)$ (according to the
result of Hong and Yang \cite{hy}), we have $v_p(P_k) \leq
v_p(\lcm(1 , 2 , \dots , k))$, that is $v_p(P_k) \leq \max_{1 \leq
i \leq k} v_p(i)$. So there exists $i_0 \in \{1 , 2 , \dots , k\}$
such that $v_p(P_k) \leq v_p(i_0)$. It follows, according to the
elementary properties of the $p$-adic valuation, that we have:
$$v_p(P_k) = \min\left(v_p(P_k) , v_p(i_0)\right) \leq v_p(P_k + i_0) \leq \max_{1 \leq i \leq k} v_p(P_k + i) ,$$
which confirms (\ref{eq2}) and completes this
proof.\penalty-20\null\hfill$\blacksquare$\par\medbreak

\noindent{\bf Proof of Theorem \ref{t7}.} Let $k \in \mathbb{N}$
fixed. The main idea of the proof is to calculate in two different ways the
quotient $\frac{g_k(P_k)}{g_k(P_k + 1)}$ and then to compare the
obtained results. On one hand, we have from the definition of the
function $g_k$:
\begin{eqnarray}\label{eq3}
\frac{g_k(P_k)}{g_k(P_k + 1)} & \!\!\!\!=\!\!\!\! & \frac{P_k (P_k
+ 1) \dots (P_k + k)}{\lcm(P_k , P_k + 1 , \dots , P_k + k)} /
\frac{(P_k + 1) (P_k + 2) \dots (P_k + k + 1)}{\lcm(P_k + 1 , P_k
+ 2 , \dots , P_k + k + 1)} \notag \\ & \!\!\!\!=\!\!\!\! & P_k
\frac{\lcm(P_k + 1 , P_k + 2 , \dots , P_k + k + 1)}{(P_k + k + 1)
\lcm(P_k , P_k + 1 , \dots , P_k + k)}
\end{eqnarray}
Next, using Lemma \ref{l1} and the well-known formula ``$a b =
\lcm(a , b) \gcd(a , b)$ $(\forall a , b \in \mathbb{N}^*)$'', we
have:
$$(P_k + k + 1) \lcm(P_k , P_k + 1 , \dots , P_k + k) = (P_k + k + 1) \lcm(P_k + 1 , P_k + 2 , \dots , P_k + k)$$
 $$= \lcm\left(P_k + k + 1 , \lcm(P_k + 1 , \dots , P_k + k)\right)~~~~~~~~~~~~~~~~~~~~~~~~~~~~~~~~~~~~~~~~~~~~
 ~~~~~~~~~~~~~~~~~~~~~~~~~~~~~~~~$$
 $$\times \gcd\left(P_k + k + 1 , \lcm(P_k + 1 ,
\dots , P_k + k)\right) ~~~~~~~~~~~~~~~~~~~~~~~~~~~~~~~~~~~~~~~~$$
 $$= \lcm(P_k + 1 , P_k + 2 , \dots , P_k + k + 1) \gcd\left(P_k + k + 1 , \lcm(P_k + 1 , \dots , P_k + k)\right)
 .$$
By substituting this into (\ref{eq3}), we obtain:
\begin{equation}\label{eq4}
\frac{g_k(P_k)}{g_k(P_k + 1)} = \frac{P_k}{\gcd\left(P_k + k + 1 ,
\lcm(P_k + 1 , \dots , P_k + k)\right)} .
\end{equation}
On other hand, according to Proposition \ref{p} and to the
definition of $P_k$, we have:
\begin{equation}\label{eq5}
\frac{g_k(P_k)}{g_k(P_k + 1)} = \frac{k!}{g_k(1)} = \frac{\lcm(1 ,
2 , \dots , k + 1)}{k + 1} .
\end{equation}
Finally, by comparing (\ref{eq4}) and (\ref{eq5}), we get:
$$P_k = \frac{\lcm(1 , 2 , \dots , k + 1)}{k + 1} \gcd\left(P_k + k + 1 , \lcm(P_k + 1 , P_k + 2 ,
\dots , P_k + k)\right) ,$$ as required. The proof is
complete.\penalty-20\null\hfill$\blacksquare$\par\medbreak

From Theorem \ref{t7}, we derive the following interesting
corollary, which confirms the conjecture of Hong and Yang
\cite{hy}.
\begin{coll}\label{c2}
For all $k \in \mathbb{N}$, the exact period $P_k$ of $g_k$ is a
multiple of the positive integer $\frac{\lcm(1 , 2 , \dots , k , k
+ 1)}{k + 1}$. In addition, for all $k \in \mathbb{N}$ for which
$(k + 1)$ is prime, we have precisely $P_k = \lcm(1 , 2 , \dots ,
k)$.
\end{coll}
{\bf Proof.} The first part of the corollary immediately follows
from Theorem \ref{t7}. Furthermore, we remark that if $k$ is a
natural number such that $(k + 1)$ is prime, then we have
$\frac{\lcm(1 , 2 , \dots , k + 1)}{k + 1} = \lcm(1 , 2 , \dots ,
k)$. So, $P_k$ is both a multiple and a divisor of $\lcm(1 , 2 ,
\dots , k)$. Hence $P_k = \lcm(1 , 2 , \dots , k)$. This finishes
the proof of the
corollary.\penalty-20\null\hfill$\blacksquare$\par\medbreak

Now, we exploit the identity of Theorem \ref{t7} in order to
obtain the $p$-adic valuation of $P_k$ $(k \in \mathbb{N})$ for
most prime numbers $p$.
\begin{thm}\label{t8}
Let $k \geq 2$ be an integer and $p \in [1 , k]$ be a prime number
satisfying:
\begin{equation}\label{eq6}
v_p(k + 1) < \max_{1 \leq i \leq k} v_p(i) .
\end{equation}
Then, we have:
$$v_p(P_k) = \max_{1 \leq i \leq k} v_p(i) .$$
\end{thm}
{\bf Proof.} The identity of Theorem \ref{t7} implies the
following equality:
\begin{equation}\label{eq7}
v_p(P_k) = \max_{1 \leq i \leq k + 1} (v_p(i)) - v_p(k + 1) +
\min\left\{v_p(P_k + k + 1) , \max_{1 \leq i \leq k}(v_p(P_k +
i))\right\} .
\end{equation}
Now, using the hypothesis (\ref{eq6}) of the theorem, we have:
\begin{equation}\label{eq8}
\max_{1 \leq i \leq k + 1}(v_p(i)) = \max_{1 \leq i \leq
k}(v_p(i))
\end{equation}
and
\begin{equation*}
\max_{1 \leq i \leq k + 1}(v_p(i)) - v_p(k + 1) > 0 .
\end{equation*}
According to (\ref{eq7}), this last inequality implies that:
\begin{equation}\label{eq9}
\min\left\{v_p(P_k + k + 1) , \max_{1 \leq i \leq k} v_p(P_k +
i)\right\} < v_p(P_k) .
\end{equation}
Let $i_0 \in \{1 , 2 , \dots , k\}$ such that $\max_{1 \leq i \leq
k} v_p(i) = v_p(i_0)$. Since $P_k$ divides $\lcm(1 , 2 , \dots ,
k)$, we have $v_p(P_k) \leq v_p(i_0)$, which implies that $v_p(P_k
+ i_0) \geq \min(v_p(P_k) , v_p(i_0)) = v_p(P_k)$. Thus $\max_{1
\leq i \leq k} v_p(P_k + i) \geq v_p(P_k)$. It follows from
(\ref{eq9}) that
\begin{equation}\label{eq10}
\min\left\{v_p(P_k + k + 1) , \max_{1 \leq i \leq k} v_p(P_k +
i)\right\} = v_p(P_k + k + 1) < v_p(P_k) .
\end{equation}
So, we have
$$\min\left(v_p(P_k) , v_p(k + 1)\right) \leq v_p(P_k + k + 1) < v_p(P_k) ,$$
which implies that
$$v_p(k + 1) < v_p(P_k)$$
and then, that
$$v_p(P_k + k + 1) = \min\left(v_p(P_k) , v_p(k + 1)\right) = v_p(k + 1) .$$
According to (\ref{eq10}), it follows that
\begin{equation}\label{eq11}
\min\left\{v_p(P_k + k + 1) , \max_{1 \leq i \leq k} v_p(P_k +
i)\right\} = v_p(k + 1) .
\end{equation}
By substituting (\ref{eq8}) and (\ref{eq11}) into (\ref{eq7}), we
finally get:
$$v_p(P_k) = \max_{1 \leq i \leq k} v_p(i) ,$$
as required. The theorem is
proved.\penalty-20\null\hfill$\blacksquare$\par\medbreak

Using Theorem \ref{t8}, we can find infinitely many natural
numbers $k$ so that $(k + 1)$ is not prime and the equality $P_k =
\lcm(1 , 2 , \dots , k)$ holds. The following corollary gives
concrete examples for such numbers $k$.
\begin{coll}\label{c4}
If $k$ is an integer having the form $k = 6^r - 1$ $(r \in
\mathbb{N}, r \geq 2)$, then we have
$$P_k = \lcm(1 , 2 , \dots , k) .$$
Consequently, there are an infinitely many $k \in \mathbb{N}$ for
which $(k + 1)$ is not prime and the equality $P_k = \lcm(1 , 2 ,
\dots , k)$ holds.
\end{coll}
{\bf Proof.} Let $r \geq 2$ be an integer and $k = 6^r - 1$. We
have $v_2(k + 1) = v_2(6^r) = r$ while $\max_{1 \leq i \leq k}
v_2(i) \geq r + 1$ (since $k \geq 2^{r + 1}$). Thus $v_2(k + 1) <
\max_{1 \leq i \leq k} v_2(i)$.\\
Similarly, we have $v_3(k + 1) = v_3(6^r) = r$ while $\max_{1 \leq
i \leq k} v_3(i) \geq r + 1$ (since $k \geq 3^{r + 1}$). Thus
$v_3(k + 1) < \max_{1 \leq i \leq k} v_3(i)$.\\
Finally, for any prime $p \in [5 , k]$, we clearly have $v_p(k +
1) = v_p(6^r) = 0$ and $\max_{1 \leq i \leq k} v_p(i) \geq 1$.
Hence $v_p(k + 1) < \max_{1 \leq i \leq k} v_p(i)$.\\
This shows that the hypothesis of Theorem \ref{t8} is satisfied
for any prime number $p$. Consequently, we have for any prime $p$:
$v_p(P_k) = \max_{1 \leq i \leq k} v_p(i) = v_p(\lcm(1 , 2 , \dots
, k))$. Hence $P_k = \lcm(1 , 2 , \dots , k)$, as
required.\penalty-20\null\hfill$\blacksquare$\par\medbreak

\section{Determination of the exact value of $P_k$}~

Notice that Theorem \ref{t8} successfully computes the value of
$v_p\left(P_k\right)$ for almost all primes $p$ (in fact we will
prove in Proposition \ref{nbadprimes} that Theorem \ref{t8} fails
to provide this value for at most one prime).  In order to
evaluate $P_k$, all we have left to do is compute
$v_p\left(P_k\right)$ for primes $p$ so that $v_p(k+1)\geq
\max_{1\leq i\leq k} v_p(i)$.  In particular we will prove:

\begin{lemma}\label{divprimeslem}
Let $k \in \mathbb{N}$. If $v_p(k+1)\geq \max_{1\leq i\leq k}
v_p(i)$, then $v_p\left( P_k\right) = 0.$
\end{lemma}

From which the following result is immediate:
\begin{thm}\label{mainthm} We have for all $k \in \mathbb{N}$:
$$
P_k = \prod_{p \ \textrm{prime}, \ p\leq k} p^{\begin{cases}0 &
\textrm{if} \ v_p(k+1)\geq \max_{1\leq i\leq k} v_p(i)\\
\max_{1\leq i\leq k} v_p(i) & \textrm{else} \end{cases}}.
$$
\end{thm}

In order to prove this result, we will need to look into some of
the more detailed divisibility properties of $g_k(n)$.  In this
spirit we make the following definitions:

Let $S_{n,k}=\{n,n+1,n+2,\ldots,n+k\}$ be the set of integers in the range $[n,n+k]$.

For a prime number $p$, let $g_{p,k}(n) := v_p(g_k(n))$. Let
$P_{p,k}$ be the exact period of $g_{p,k}$. Since a positive
integer is uniquely determined by the number of times each prime
divides it, $P_k = \lcm_{p \ \textrm{prime}}(P_{p,k}).$

Now note that
\begin{align*}
g_{p,k}(n) & = \sum_{m\in S_{n,k}} v_p(m) - \max_{m\in S_{n,k}} v_p(m)\\
& = \sum_{e>0,m\in S_{n,k}} (1 \ \textrm{if} \ p^e|m) - \sum_{e>0} (1 \ \textrm{if}
\ p^e \ \textrm{divides some} \ m\in S_{n,k})\\
& = \sum_{e>0} \max(0,\# \{m\in S_{n,k}: p^e|m\}-1).
\end{align*}
Let $e_{p,k} = \floor{\log_p (k)} = \max_{1\leq i \leq k}v_p(i)$
be the largest exponent of a power of $p$ that is at most $k$.
Clearly there is at most one element of $S_{n,k}$ divisible by
$p^e$ if $e>e_{p,k}$, therefore terms in the above sum with
$e>e_{p,k}$ are all 0.  Furthermore, for each $e\leq e_{p,k}$, at
least one element of $S_{p,k}$ is divisible by $p^e$.  Hence we
have that
\begin{equation}\label{gpkformula}
g_{p,k}(n) = \sum_{e=1}^{e_{p,k}} \left(\# \{m\in S_{n,k}:
p^e|m\}-1\right).
\end{equation}

Note that each term on the right hand side of (\ref{gpkformula})
is periodic in $n$ with period $p^{e_{p,k}}$ since the condition
$p^e|(n+m)$ for fixed $m$ is periodic with period $p^e$.
Therefore $P_{p,k}|p^{e_{p,k}}$.  Note that this implies that the
$P_{p,k}$ for different $p$ are relatively prime, and hence we
have that
$$
P_k = \prod_{p \ \textrm{prime}, \ p\leq k} P_{p,k}.
$$

We are now prepared to prove our main result

\noindent{\bf Proof of Lemma \ref{divprimeslem}.} Suppose that
$v_p(k+1) \geq e_{p,k}$.  It clearly suffices to show that
$v_p\left(P_{q,k}\right)=0$ for each prime $q$.  For $q\neq p$
this follows immediately from the result that
$P_{q,k}|q^{e_{q,k}}$.  Now we consider the case $q=p$.

For each $e \in \{1 , \dots , e_{p , k}\}$, since $p^e | k + 1$,
it is clear that $\# \{m\in S_{n,k}: p^e|m\} = \frac{k+1}{p^e}$,
which implies (according to (\ref{gpkformula})) that $g_{k , n}$
is independent of $n$. Consequently, we have $P_{p,k}=1$, and
hence $v_p\left( P_{p,k}\right) = 0$, thus completing our proof.
\penalty-20\null\hfill$\blacksquare$\par\medbreak Note that a
slightly more complicated argument allows one to use this
technique to provide an alternate proof of Theorem \ref{t8}.

We can also show that the result in Theorem \ref{mainthm} says that $P_k$ is basically $\lcm(1,2,\ldots,k)$.
\begin{prop}\label{nbadprimes}
There is at most one prime $p$ so that $v_p(k+1)\geq e_{p,k}$. In
particular, by Theorem \ref{mainthm}, $P_k$ is either
$\lcm(1,2,\ldots,k)$, or $\frac{\lcm(1,2,\ldots,k)}{p^{e_{p,k}}}$
for some prime $p$.
\end{prop}
\noindent {\bf Proof.} Suppose that for two distinct primes, $p,q \leq k$ that $v_p(k+1)\geq e_{p,k}$, and $v_q(k+1)\geq e_{q,k}$.
Then
$$
k+1 \geq p^{v_p(k+1)}q^{v_q(k+1)} \geq p^{e_{p,k}}q^{e_{q,k}} > \min\left(p^{e_{p,k}},q^{e_{q,k}} \right)^2 = \min\left(p^{2e_{p,k}},q^{2e_{q,k}} \right).
$$
But this would imply that either $k\geq p^{2e_{p,k}}$ or that $k\geq q^{2e_{q,k}}$ thus violating the definition of either $e_{p,k}$ or $e_{q,k}$.
\penalty-20\null\hfill$\blacksquare$\par\medbreak

\end{document}